\documentclass[11pt]{article}
\usepackage[dvips]{graphicx}
\usepackage{a4wide}
\usepackage{wasysym}
\usepackage{amsfonts}
\usepackage{amssymb}
\usepackage{epsfig}

\usepackage[all]{xy}
\newtheorem{defi}{\textbf{Definition}}[section]
\newtheorem{theo}[defi]{\textbf{Theorem}}
\newtheorem{lemma}[defi]{\textbf{Lemma}}
\newtheorem{prop}[defi]{\textbf{Proposition}}
\newtheorem{coro}[defi]{\textbf{Corollary}}

\newcommand{\Diff}{{\rm Diff}}
\newcommand{\Diam}{{\rm Diam}}
\newcommand{\length}{{\it length}}
\newcommand{\Card}{\sharp}
\newcommand{\MCG}{\mathcal{MCG}}
\newcommand{\CG}{\mathcal{CG}}
\newcommand{\cptg}{\mathfrak{cptg}}
\newcommand{\Capa}{\kappa}

\title{Copies of a one-ended group in a Mapping Class Group}

\author{Fran\c{c}ois  Dahmani \footnote{The first author acknowledges partial support from the ANR grant ANR-06-JCJC-0099 }, 
Koji Fujiwara 
\footnote{The second author is partially supported
by  Grant-in-Aid for Scientific Research (No. 19340013).}
}
\date{}
\begin{document} 

\maketitle

\begin{center}{\bf Abstract. } We establish that, given $\Sigma$ a compact orientable surface, and $G$ a finitely presented one-ended group, the set of copies of $G$ in the mapping class group $\mathcal{MCG}(\Sigma)$ consisting of only pseudo-Anosov elements except identity, is finite up to conjugacy. This relies on a result of Bowditch  on the same problem for images of surfaces groups.  He asked us whether we could reduce the case of one-ended groups to his result ; this is a positive answer. Our work involves analogues of Rips and Sela's canonical cylinders in curve complexes, and an argument of Delzant to bound the number of images of a group in a hyperbolic group.
\end{center}

 Let $\Sigma$ be a compact orientable surface (possibly with boundary components). 
The Mapping Class Group of $\Sigma$, denoted by
$\mathcal{MCG}(\Sigma)$ is the group of isotopy classes of orientation
preserving  self-homeomorphisms of $\Sigma$.

 The aim of this paper is to report on a control on the family of subgroups of
Mapping Class Groups that are isomorphic to a given finitely presented one-ended group.

\begin{defi}[Purely pseudo-Anosov]
 A subgroup of $\mathcal{MCG}(\Sigma)$ is said purely pseudo-Anosov if all its elements, except the identity, are pseudo-Anosov mapping classes. A morphism of a group in  $\mathcal{MCG}(\Sigma)$ is said purely pseudo-Anosov if it is injective and has  purely pseudo-Anosov image.
\end{defi}

Up to now, the only known purely pseudo-Anosov subgroups of Mapping Class Groups are free. 

Recently, Brian Bowditch \cite{B3} has established the finiteness of the set of such images of surface groups, up to conjugacy.  He uses deep results, in particular in 3-manifold geometry, and doing so, he proves a general powerful result   \cite[Proposition 8.1]{B3}.  He asked us, however, whether it is possible to adapt the situation of an arbitrary finitely presented one-ended group to the setting of his Proposition. We provide here an affirmative answer.

\begin{theo}\label{theo;intro1}
  Given $\Sigma$ a compact orientable surface, and $G$ a finitely presented one-ended group, there is a number
  $N$ such that any purely pseudo-Anosov subgroup of $\mathcal{MCG}(\Sigma)$ isomorphic to $G$
  admits a set of generators  $(\gamma_i)$ for which there are a vertex $v$ in the curve complex $\mathcal{CC}(\Sigma)$
  with $d(\gamma_i v, v) \leq N$, and a presentation on this set of generators with at most $N$ relators, each of them of length at most $N$ as words.
\end{theo}

This allows to apply the following important Bowditch's result: 
\begin{prop} \cite[Proposition 8.1]{B3} 
 Suppose that $G$ is a one-ended finitely presented group, and that
  $\phi: G \to  \mathcal{MCG}(\Sigma)$ is a purely pseudo-Anosov  homomorphism, giving an induced action
on $\mathcal{CC}(\Sigma)$. Suppose that  $A$   is a generating set
of $G$  and that there is some $v\in \mathcal{CC}(\Sigma)$ and    $N>0$
 such that $d(\phi(a)v,v) \leq N$  for all $a \in A$. Then there is some  $\theta \in  \mathcal{MCG}(\Sigma)$ such that the word
length of each $ \theta \phi(a) \theta^{-1}, a\in A$ 
 (in terms of a generating set of 
  $\mathcal{MCG}(\Sigma)$) is bounded above in terms of $N$ and the sum of lengths of relators in a presentation of $G$ on the generating set $A$. 
\end{prop}

The following corollary is immediate
since for each given $G$ and $A$, the set of 
the elements $ \theta \phi(a) \theta^{-1}, a \in A$
for all $\phi$ is finite if we choose $\theta$
using the proposition for each $\phi$.

\begin{coro}
 Given $\Sigma$ a compact orientable surface, and $G$ a finitely presented one-ended group, the set of purely pseudo-Anosov subgroups of $\mathcal{MCG}(\Sigma)$ isomorphic to $G$ is finite up to conjugacy in $\mathcal{MCG}(\Sigma)$.
\end{coro}

 The group $\mathcal{MCG}(\Sigma)$ has a natural action  by isometries,  on Harvey's
 curve complex $\mathcal{CC}(\Sigma)$,
which turns out to be a hyperbolic space \cite{MM}, \cite{B1}. This complex is far from being
locally finite, and the action is not proper. 

The elements of $\mathcal{MCG}(\Sigma)$ that are hyperbolic isometries of   $\mathcal{CC}(\Sigma)$    are precisely the  pseudo-Anosov elements of $\mathcal{MCG}(\Sigma)$.

Our method toward Theorem \ref{theo;intro1} is inspired by the case of relative hyperbolicity, studied in \cite{Dis} (and indeed to the hyperbolic case, \cite{Del}): construct Rips and Sela's canonical
cylinders, as in  \cite{RS}, for a given morphism $G \to \mathcal{MCG}(\Sigma)$, and
use them to pull back a lamination on a Van Kampen complex $P(G)$ of $G$ (or more precisely, first on its universal cover), that allows to find small
generators of the image.

To perform the construction, we make use of 
 deep results of  Bowditch about tight geodesics in $\mathcal{CC}(\Sigma)$.

In this  paper, we introduce the relevant definitions for the general argument, but sometimes refer to existing proofs when they can be applied without modification. We tried to make clear where the technology specific to Mapping Class Groups is used, or where the existing argument would not, as written, give sufficient precision.
In particular, our main task about Theorem \ref{theo;canonicalcylinders}, which is based on very subtle ideas of Rips and Sela, is to explain how to get a setting where the original proof can be applied (which is not obvious without Bowditch's results). However, for the reader's convenience, we also reproduce this proof in section \ref{sec;proof1.9}.

Let us mention that in the case  $G$ is also a surface group, other proofs of Theorem \ref{theo;intro1} have been given, notably by J.~Barnard \cite{Ba}.

We thank B.~Bowditch for stimulating discussions, and for asking the question on the bound on the complexity of presentations. We learned, while finishing this paper,  that he very simultaneously obtained a similar result, by different methods, using actions on $\mathbb{R}$-trees \cite{B4}.  The first author thanks T.~Delzant  for interesting and stimulating related discussions. 
Both authors want to thank the referee for constructive remarks.
 The second author gratefully acknowledges the
Institut de Math\'ematiques de Toulouse
and support from CNRS. 
This work was partially 
done while he was visiting the institute.

\section{Sliced canonical cylinders}

 In the following, $\mathcal{CG}(\Sigma)$ is the curve graph of a surface $\Sigma$ (the one-skeleton of Harvey's curve complex $\mathcal{CC}(\Sigma)$),  
 $\delta\in \mathbb{N}$ is an hyperbolicity constant and
$p$ is a base point in $\mathcal{CG}(\Sigma)$.

 The graph   $\mathcal{CG}(\Sigma)$ is not locally finite, and in general two points are joined by infinitely many different geodesics.  However, there is a class of geodesics that are called \emph{tight geodesics}, and that have good properties.  We will not need the definition, which involves properties of the curves  and subsurfaces in the surface $\Sigma$, so  we just refer to \cite{B2} for it.  We will need the fact that there exist such geodesics, that they satisfy the statement of Theorem \ref{theo;bowditch}, and that  a sub-path of a tight geodesic is a tight geodesic.

\begin{defi}[$\lambda$-quasi-geodesic]
A $\lambda$-quasi-geodesic in a graph $X$ is a $\lambda$-bi-lipschitz embedding of a segment of $\mathbb{R}$ into $X$. We assume here that paths start and end on  vertices. The length of a path  is the number of edges in its image.

A path is a  $\nu$-local-$\lambda$-quasi-geodesic  if
 any of its subpaths of length at most $\nu$  
 is a $\lambda$-quasi-geodesic. 

A path  is a $\mu$-local tight geodesic if any of its subpaths of length
$\mu$, is a tight geodesic.

\end{defi}

We choose some constants:  $\lambda = 1000 \delta$, and $\epsilon$ such that any
$\lambda$-quasi geodesic in $\mathcal{CG}(\Sigma)$ stays $\epsilon$-close to any 
geodesic joining its end points. Let $\mu =(100 \epsilon +  \lambda^2)\times 40 \lambda$, and $\nu=40\lambda(\epsilon+100\lambda\delta)$.

The next definition is to be compared with a similar one in \cite{RS}, for geodesics that are not necessarily tight.

\begin{defi}[Coarse piecewise tight geodesics, or $\mathfrak{cptg}$]\label{def;cpg}
 Let $l \geq \mu$ be an integer.
  An $l$-coarse piecewise tight geodesic, or $l$-$\mathfrak{cptg}$, in $\mathcal{CG}(\Sigma)$ is a
$\nu$-local $\frac{\lambda}{2}$-quasigeodesic $f:[a,b] \to \mathcal{CG}(\Sigma)$ together
with a subdivision of the segment $[a,b]$, 
$a=c_1\leq d_1 \leq c_2 \dots \leq d_n =b$ such that, for all $i \leq n$,  
 $f([c_i, d_i])$ is a $\mu$-local tight geodesic, of length at least $l$ when
  $2 \leq  i \leq (n-1)$, and such that  with $\length(f[d_i, c_{i+1}])\leq \epsilon$. 

Moreover we require that $f([a,b])$ is in the $2\epsilon$-neighborhood of a tight geodesic segment $[f(a),f(b)]$.

\end{defi}

\emph{Remark: }  Any $l$-$\mathfrak{cptg}$  is a $\lambda$-quasi-geodesic
 (this does not use tightness, only hyperbolic geometry; see for instance \cite[Appendix]{Dis}). The subpaths corresponding to a subdivision of an $l$-$\mathfrak{cptg}$  as in the definition are called respectively sub-local geodesics, and bridges.

\begin{defi}[$l$-Cylinders, \cite{RS}]\label{def;cylinders}
Let $l\in \mathbb{N}$.
 The $l$-cylinder of two points $x$ and  $y$ in
$\mathcal{CG}(\Sigma)$, denoted by $Cyl_l(x, y)$, 
is the set of the vertices $v$ lying on
an $l$-$\mathfrak{cptg}$ from $x$ to $y$, with  the additional
requirement that $v$ is on a local tight geodesic $f|_{[c,d]}$ of the subdivision, with distances
$|f(c) - v|\geq l$ if $f(c)\neq x$  and $|f(d)-v| \geq l$ if $f(d)\neq y$.
\end{defi}

  Any tight geodesic is an $l$-$\mathfrak{cptg}$, 
  with a trivial subdivision, and thus is contained in the $l$-cylinder 
  between its end points, for all $l$. 
Here is an obvious consequence of definitions.

\begin{lemma}[Equivariance]
 $\forall \gamma \in \MCG(\Sigma),\,  \forall x,y, \;  
\gamma Cyl_l(x,y) =  Cyl_l(\gamma x, \gamma y)$. Moreover $Cyl_l(x,y)= Cyl_l(y,x)$
\end{lemma}

Recall a crucial result of Bowditch (that will be used in the next two lemmas):
\begin{theo}\label{theo;bowditch}

\cite[Theorem 1.1]{B2}: Given $L$, there is a constant $K_0(L)$ such that for all $a,b$ vertices of the curve complex, and all $c$ on a tight geodesic between $a$ and $b$,  the set of vertices on tight geodesics between $a$ and $b$ and at distance at most $ L+2\epsilon$ from $c$, has at most $K_0(L)$ elements.

\cite[Theorem 1.2]{B2}:
There are constants $k_1$ and $K_1$ depending only on $\Sigma$ such that if $a,b$ are vertices in $\mathcal{CC}(\Sigma)$,  $r\in \mathbb{N}$,  and $c$ on a tight geodesic joining $a$ to $b$, with $d(a,c) \ge r + k_1$ and $d(b,c) \ge r +k_1$,  then the set of vertices on tight geodesics between two points respectively $r$-close of $a$ and $b$, and at distance at most $2\epsilon$ from $c$, has at most $K_1$ elements. 

\end{theo}

\begin{lemma}\label{lem;cylfinite}
Let $l\geq  k_1+\mu/2  $.  Any $l$-cylinder of the curve complex  $ \mathcal{CC}(\Sigma)$ is finite.
\end{lemma}

{\it Proof. } Let $x,y$ be two vertices in the curve complex. Let $C$ be the set of vertices  on tight geodesics between $x,y$ that are at distance at least $l- \mu/4$ from both $x$ and $y$. By the first point of Theorem \ref{theo;bowditch}, this set is finite. For each $v\in C$ let $B_v$ be the set given by the second point of Theorem \ref{theo;bowditch} for $r=\mu/4$ (which can be applied since $l$ is large enough). Let $B$ be the union of all the $B_v$, $v\in C$, this set is finite.  Let also $B'$ be the set of all vertices on tight geodesics from $x$ to $y$ or from $x$ to a vertex of $B$, or from $y$ to a vertex of $B$. Since $B$ is finite, and by the first point of Theorem \ref{theo;bowditch}, this set $B'$ is finite. We want to show that the cylinder of $x,y$ is a subset of $B\cup B'$.

Let $w$ lie in an $l$-cylinder of $x,y$.  It is on a local tight geodesic $f|_{[c,d]}$ with $f(c)$ and $f(d)$ at distance at most $2\epsilon$ from a tight geodesic between $x$ and $y$.  

Thus, $w$ is on a subsegment $\sigma$ of length $\mu$ that is a tight geodesic, and whose end points, are at distance at most $2\epsilon$ from a tight geodesic $[x,y]$. There are two cases following from the inequality condition of Definition \ref{def;cylinders}: either we can assume that one of the ends of $\sigma$ is $x$ (or $y$) and $d(w,x) \leq \mu/2$ (or similarly with $y$), or we can assume that  $w$ is in the middle of $\sigma$.

 In the second case, since $\mu >100\epsilon$, one can find a smaller   sub-(tight geodesic) of length $\mu/2$ with  $w$ in its middle, and whose end points are, by hyperbolicity, at distance at most $2\delta (<\mu/4)$ from  a tight geodesic $[x,y]$.  Then by definition of the sets $B_v$, we have that  $w$ in some set $B_v$ for some $v\in C$.

In the first case, assume that one end of $\sigma$  is $x$. Then applying the above argument to the center $c(\sigma)$ of $\sigma$, we find that  $c(\sigma) \in B$, and therefore, $w\in B'$. $\square$

\begin{defi}[Channels, compare to \cite{RS}, 4.1]
Let $L>0$, and   $a,b \in \mathcal{CG}(\Sigma)$ with $d(a,b) \leq 3L$.  A tight geodesic $g_1$ of length $L$ which is contained in a tight geodesic $g_2$ of length $3L$ that  starts (respectively  ends) at distance at most $2\epsilon$ from  $a$ (respectively $b$), such that end points of  $g_2$ are at distance $L$ from  $g_1$,   
is called an $L$-\emph{channel}
of $(a,b)$. \end{defi}

\begin{lemma}\label{lem;chan}
 For every $L\ge k_1 +4\epsilon$,   
there is a bound $\Capa(L)$ on the number of  $L$-channels of $(a,b)$ for arbitrary $a,b$ with $d(a,b) \leq 3L$.

\end{lemma}

{\it Proof. } First we can assume that $d(a,b) \geq 3L -4\epsilon$ otherwise there are no channels at all. 

We will show a finite set (of cardinality uniformly bounded above) containing all vertices of $L$-channels of $a,b$.
Let $C$ be the set of vertices on tight geodesics from $a$ to $b$ that are at distance at least $L-2\epsilon$ from both. Because $d(a,b) \leq 3L$, by the first point of Theorem \ref{theo;bowditch}, $C$ has at most $K_0(L)$ elements.
For each $c\in C$, consider $B_c$ the set given by the second point of Bowditch's theorem, for $r=2\epsilon$,  applied to $c$ (which can be applied  since $d(a,c) \geq L-2\epsilon \geq k_1 +2\epsilon$ by choice of $L$). Let $B$ be the union of all the $B_c ,c\in C$. It has at most $K_0\times K_1$ elements.

Now consider $w$ a vertex on an $L$-channel $g_1$ of $a,b$, which by definition is a sub-(tight geodesic) of a tight geodesic $g_2$ starting and ending at $a',b'$ with $d(a,a' ) \leq 2\epsilon$ and $d(b,b' ) \leq 2\epsilon$.  Let $L+x = d(a',w)$ $(x\geq 0)$, and assume (by symmetry this is without loss of generality) that $x\leq L/2$.  We have $   L-2\epsilon +x \leq d(a,w) \leq L+2\epsilon +x$. By hyperbolicity, $w$ is $2 \delta$-close to a point $w'$ in a segment $[a,b]$ (which we choose tight). Since $2\delta \leq \epsilon$, we can find another point $w''$ on $[a,b]$ at distance at most $\epsilon$ from $w$  (hence $w\in B_{w''}$)  and at least $  L-2\epsilon +x $ from $a$ (hence $w''\in C$). This shows that $w\in B$.  $\square$

For an integer $n$, we set
$\psi(n) = 24(n+1)\Capa(\mu)(2\epsilon+1)\epsilon $.
We denote by $B_r(x)$ the ball of $\mathcal{CG}(\Sigma)$ of center $x$ and radius $r$.

\begin{theo}\label{theo;canonicalcylinders}

 Let $F$ be a finite family of elements of $\MCG$ ; we set $n = (2
\Card(F))^3$ where $\Card(F)$ is the cardinality of $F$.
 Let $p$ be a base point
in $\mathcal{CG}(\Sigma)$.

  There exists a number $l$ such that the $l$-cylinders satisfy: for all
$\alpha, \beta, \gamma $ in  $F\cup F^{-1}$ with $\alpha\beta\gamma =1$, in the
triangle $(x,y,z)=(p,\alpha p, \gamma^{-1} p)$ in $\mathcal{CG}(\Sigma)$, one has 
$$ 
Cyl_l(x,y) \cap B_{R_{x,y,z}}(x) =  Cyl_l(x,z) \cap B_{R_{x,y,z}}(x)
$$
 (and
analogues permuting $x$, $y$ and $z$) where $R_{x,y,z} = (y\cdot z)_x -
5\times (13 \mu+\psi(n))$, is the Gromov product in the triangle,
minus a constant.

\end{theo}

Note that in the statement, $l$ depends on $F$, but $(y\cdot z)_x - R_{x,y,z}$  depends only on  $\Card(F)$.

 One may think that $Cyl_l(x,y)$ is a narrow
set near a geodesic from $x$ to $y$.
The theorem says that $Cyl_l(x,y), Cyl_l(y,z)$
and $Cyl_l(z,x)$ coincide except in a set of bounded 
size near the center of the triangle $(x,y,z)$.
Instead of $Cyl_l(x,y)$, if we take the union 
of all (tight) geodesics between $x,y$ or 
the union of all quasi-geodesics between $x,y$ with 
uniform quasi-geodesic constants, we do not  
have this equation in general. This is already the case for 
a Cayley graph of a word-hyperbolic group, 
and Rips-Sela \cite{RS} introduced several notions in this context 
which we imitated here.

\subsection{Proof of Theorem \ref{theo;canonicalcylinders}}\label{sec;proof1.9}

We produced a setting where cylinders and channels are finite. We can therefore reproduce the original proof of Theorem  \ref{theo;canonicalcylinders} by Rips and Sela for hyperbolic groups \cite{RS}. For the reader's convenience, we give the detail. We follow the exposition in  \cite[Theorem 2.9]{Dis} (which was for relatively hyperbolic groups).

Let us start by stating two lemmas for rerouting a $\cptg$.

\begin{lemma}\label{lem;reroute1}
 Let $l\geq \mu$, and $f:[a,b]\to \CG$ be an $l$-$\cptg$, whose subdivision includes $f|_{[c,d]}$, a local tight geodesic. Let $s=f(t)\in f([c,d])$ such that the path  $f([c,t])$ from $f(c)$ to $s$ has length $\geq l+2\epsilon$.

Let now $g$ be a tight geodesic segment joining $f(a)$ to  $f(b)$ and $s''$ be a point on $g$ closest to $s$. Let $s' = f(t')\in f( [c,d])$  be a point closest to $s''$ on 
$f([c,d])$. 

Let $[s',s'']$ be a geodesic segment, and $[s'',f(b)]$ be a subsegment of $g$.  

Then, the concatenation $f\left([a,t']\right) \ast [s',s'']\ast [s'',f(b)]$ is 
an $l$-$\cptg$ from 
$f(a)$ to $f(b)$, and we say that $f$ can be rerouted into this new path.
  
\end{lemma}

\begin{lemma}\label{lem;reroute2}
  Let $f$ be an $l$-$\cptg$ whose last subdivision segment is a local tight geodesic $g$ of length at least $l+2\mu$. 
Let $z \in \CG$ such that   a tight geodesic segment $[f(a),z]$ passes at distance at most $\delta$ from $f(b)$. 

Then there exists an $l$-$\cptg$ from $f(a)$ to $z$ coinciding with $f$ until the first point of $g$.
\end{lemma}

\begin{figure}[hbt]
            \begin{center}
              \includegraphics{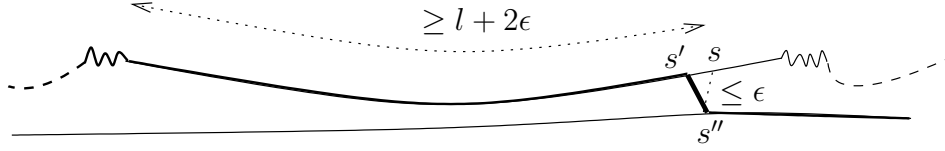}
              \caption{A typical example of rerouting a  $\cptg$ (see  Lemma \ref{lem;reroute1} and \ref{lem;reroute2}).  Starting with a  $\cptg$ one of whose  local geodesics is longer than $l+2\epsilon$,  one deduces another $\cptg$ (bold on the picture)  coinciding with the first until a certain point   on this long local geodesic ($s'$ on the picture), but ending to a possibly different point.} 
                \label{fig;reroute}
           \end{center}
         \end{figure}

We do not repeat the proofs of these two lemmas here. They are rather standard, we refer to Lemma 2.2 and 2.4 in \cite{Dis} for instance, see Figure \ref{fig;reroute} for an illustration. The main observation is that the proposed paths are indeed local quasi-geodesics (the other properties being immediate).

\begin{lemma}\label{lem;fail_to_satisfy}

For all integer $i$ in $\{1,\dots,  \frac{\psi(n)}{2\epsilon} \} $,  let  us define $l_i = 10\mu + 2i\epsilon \leq  \psi(n)+10\mu$.

 Let $x,y,z$ be three points in $\CG$. There are at most $2\Capa(\mu)\times (2\epsilon+1)$ different values of $l_i$ for  $i\in\{1, \dots ,  \frac{\psi(n)}{2\epsilon}  \}$ such that 
 $$ Cyl_{l_i}(x,y)  \cap   B_{R_{x,y,z}}(x)     \not\subset   Cyl_{l_i}(x,z) \cap  B_{R_{x,y,z}}(x).    $$ 

\end{lemma} 
 
{\it Proof. } We argue by contradiction, assuming that $2\Capa(\mu)\times (2\epsilon+1) +1$ different $l_i$ do not satisfy $ Cyl_{l_i}(x,y)  \cap   B_{R_{x,y,z}}(x)     \subset  Cyl_{l_i}(x,z) \cap  B_{R_{x,y,z}}(x)$. 

For each of them, there is $v_i \in B_{R_{x,y,z}}(x)$ in $ Cyl_{l_i}(x,y)$ but not in  $Cyl_{l_i}(x,z)$: there is $\beta_i$ an $l_i$-$\cptg$ from $x$ to $y$ containing $v_i$ as indicated in the definition \ref{def;cylinders},  and none from $x$ to $z$.

 By definition of $\cptg$,  given a tight geodesic $[x,y]$, $\beta_i$ is contained in its $2\epsilon$-neighborhood. Thus every subsegment of length $\mu$ of a sub-local geodesic of $\beta_i$, at distance at least $\mu$ from the end points of the sub-local geodesic,  is in fact a $\mu$-channel
of some subsegment of $[x,y]$. 

Let $[x',y']$ and $[x'',y'']$ be two subsegments of $[x,y]$ of length $3\mu$, such that $d(x',x) =R_{x,y,z} + (\psi(n) +9\mu)$ and $d(x'',x) =R_{x,y,z} + (\psi(n) +13\mu)$ (the end of the first of these segments  is at distance $\mu$ from the beginning of the second). Assume that $\beta_i$ 
does not contain a $\mu$-channel of   $[x',y']$, this means as we just noticed, that it  must have a bridge $\beta_i([d_j,c_{j+1}])$ at distance $3\mu+2\epsilon$ from $x'$. Since  $l_i \geq 10\mu$ (and  $\epsilon \leq \mu/10$), the sub-local geodesic after this bridge must contain  a $\mu$-channel of   $[x'',y'']$. Therefore, each $\beta_i$ contains a $\mu$-channel of either  $[x',y']$ or $[x'',y'']$

There are at most $2\Capa(\mu)$  different $\mu$-channels of either  $[x',y']$ or $[x'',y'']$, therefore there is a channel, that we denote by $Chan$, through which passes some local geodesic subdivision $\beta_i|_{[c(i),d(i)]}$ for at least $2\epsilon+2$ different indices $i$. Let $i_1 < \dots < i_{2\epsilon+2}$ be such indices.

For each $j$, let $t_{i_j} \in [c(i_j),d(i_j)]$ be the instant where $\beta_{i_j}(t_{i_j})$ exists the channel $Chan$, and $r_{i_j}$ the length of the path $\beta_{i_j}([t_{i_j},d(i_j)])$. There are three claims about the possible values of $r_{i_j}$.

{\it Claim 1:  } For any $j \in \{1, \dots, 2\epsilon +2\}$, one has $r_{i_j} < l_{i_j} +2\epsilon$. 
 
{\it Claim 2:  } If $i_j<i_k$ then $r_{i_k} <r_{i_j}$.

{\it Claim 3:  } $r_{i_1} -r_{i_{2\epsilon +2}} <2\epsilon$.

From the second claim, we deduce that all the $r_{i_j}$ are different, and from the third claim we deduce that they are integers in an interval of length $2\epsilon +1$. Since there are $2\epsilon+2$ values this is a contradiction.

We now have to prove these claims. 

For the first one, assume the contrary, and let $t_{i_j}^+ > t_{i_j}$ be a real number such that the length of   
$\beta_{i_j}([t_{i_j},t_{i_j}^+] )$ is $l_{i_j}$. 
Our assumption allows to use Lemma \ref{lem;reroute1}: one can change $\beta_{i_j}$ into another $l_{i_j}$-$\cptg$ coinciding with $[x,y]$ on a subsegment containing $[y',y]$, with $d(y',\beta_{i_j}(t_{i_j}^+ ) ) \leq 3\epsilon$. By the triangle inequality, $d(x,\beta_{i_j} (t_{i_j}) ) \leq   R_{x,y,z} + (\psi(n)+13\mu) +2\epsilon +2\mu$. 
 Therefore, $d(x,   \beta_{i_j} (t_{i_j}^+) )  \leq   R_{x,y,z} + (\psi(n)+13\mu) +2\epsilon +2\mu + (\psi(n)+10\mu) +2\epsilon$ which is 
$\leq  R_{x,y,z} + 2 (\psi(n)+13\mu)$, which is 
$\leq (z\cdot y)_x -2 (\psi(n) +13\mu)$, meaning that it is at least $l_i +2\mu$ before reaching  a point $\delta$-close to  the center of the triangle $(x,y,z)$. This allows to use  Lemma \ref{lem;reroute2}: this new $\cptg$ can be rerouted into another one, for the same $l_{i_j}$, coinciding with the beginning $\beta_{i_j}$ until after $v_{i_j}$, and ending at $z$. In particular, $v_{i_j}$ is in $Cyl_{l_{i_j}}(x,z)$, contradicting our assumption.

Now we use the first claim to prove the second. If this second claim was not true, one could change $\beta_{i_j}$ just after its passage in $Chan$ into $\beta_{i_k}$ (it is enough to notice that this new path remains an $l_{i_j}$-$\cptg$ since $i_j<i_k$). On  $\beta_{i_k}$, consider the sub-local tight geodesic of the subdivision following that of $Chan$. Because $i_k>i_j$, it is longer than $l_{i_j} +2\epsilon$;     let $\beta_{i_k} (t_{i_k}^+)$  be the point on it after travelling this distance (which is $\leq (\psi(n) +10 \mu) +2\epsilon$ in any case). As before, $d(x,\beta_{i_j} (t_{i_j}) ) \leq   R_{x,y,z} + (\psi(n)+13\mu) +2\epsilon +2\mu$. By the first claim, $r_{i_j} \leq \psi(n) +10\mu +2\epsilon$, then the next bridge is at most $\epsilon$ long, and we need to travel at most $(\psi(n) +10 \mu) +2\epsilon$ further to find  $\beta_{i_j} (t_{i_j}^+ )$. Thus, $d(x,   \beta_{i_j} (t_{i_j}^+ ))  \leq  R_{x,y,z} + (\psi(n)+13\mu) +2\epsilon +2\mu + 2(\psi(n)+10\mu) +5\epsilon$, which is $\leq  R_{x,y,z} + 3 (\psi(n) +13\mu)  $, which is   $\leq (z\cdot y)_x -2 (\psi(n) +13\mu)$. We then use, as in claim 1,  Lemma \ref{lem;reroute1} and  Lemma \ref{lem;reroute2} to obtain the same contradiction.

The third claim is again proved by contradiction: if it was false, we could change $\beta_{i_{2\epsilon +2k}}$ just after $Chan$ by substituting the remaining part of the sub-local tight geodesic of  $\beta_{i_1}$ containing $Chan$. Then, one can reroute this $\cptg$ on $[x,y]$ at distance $2\epsilon$ before the end of this sub-local geodesic, and finally, reroute it again into a $\cptg$ ending at $z$, again a contradiction.  $\square$

Now, we can prove Theorem \ref{theo;canonicalcylinders}. We need to find a good parameter $l$. We have at least $\frac{\psi(n)}{2\epsilon} = 12(n+1)\Capa(\mu) \times (2\epsilon +1)$ candidates: the parameters $l_i$  defined in Lemma \ref{lem;fail_to_satisfy}. There are at most $n$ different triangles satisfying the condition of Theorem \ref{theo;canonicalcylinders}, hence we have a system of at most $6n$ inclusions of the form  $  Cyl_{l_i}(x,y)  \cap   B_{R_{x,y,z}}(x)   \subset     Cyl_{l_i}(x,z) \cap  B_{R_{x,y,z}}(x)   $ to satisfy simultaneously. For each inclusion, by Lemma \ref{lem;fail_to_satisfy}, only $2\Capa(\mu)(2\epsilon +1)$ parameters $l_i$ fail to satisfy it. Hence, by the pigeonhole principle, one parameter satisfies all the $6n$ inclusions. $\square$

 \subsection{Slicing}

 Let us assume that $l$ satisfies the conclusion of Theorem \ref{theo;canonicalcylinders}. From now on, all cylinders will be $l$-cylinders, and we write $Cyl(a,b)$ for $Cyl_l(a,b)$.

 Let $Cyl(a,b)$ be a cylinder and $x\in
Cyl(a,b)$.
 Following \cite{RS}, we define the set  $N^{(a,b)}_R(x)$ as follows: it is the set of all
the vertices $v\in Cyl(a,b)$ such that $|a-x|<|a-v|$, and such that
 $|x-v|>100\delta$.
 Here $R$ stands for ``right'', and $N^{(a,b)}_L(x)$
is similarly defined changing the condition $|a-x|<|a-v|$ into $|a-x|>|a-v|$.
As cylinders are finite, these sets are also finite.

 Let $x,y \in Cyl(a,b)$ be in a cylinder. 
 We set $$\begin{array}{lcllll} \Diff_{a,b}(x,y) & =&  \Card(N^{(a,b)}_L(x) \setminus  N^{(a,b)}_L(y))

& - & \Card(N^{(a,b)}_L(y)  \setminus  N^{(a,b)}_L(x)) & + \\ 
  &  &  \Card(N^{(a,b)}_R(y)
\setminus  N^{(a,b)}_R(x))&   - & \Card(N^{(a,b)}_R(x) \setminus  N^{(a,b)}_R(y)) & \end{array}$$
where $\Card(X)$ is the cardinality of the set $X$. This definition makes
 sense: because of Lemma \ref{lem;cylfinite} all the sets involved are finite.

\begin{lemma} 
 $\Diff_{a,b}$ satisfies a cocycle relation: for arbitrary $x,y, z\in Cyl(a,b)$, one has  $\Diff_{a,b}(x,z) = \Diff_{a,b}(x,y) + \Diff_{a,b}(y,z)$.

 In particular, the relation $(\Diff_{a,b}(x,y)=0)$ is an equivalence relation on $Cyl(a,b)$. 

Let us say that an equivalence class  for this relation
$(\Diff_{a,b}(x,y)=0)$ is a \emph{slice} of  $Cyl(a,b)$.

 The value of  $\Diff_{a,b}(x,y)$  depends only on the slices of $x$ and $y$. 

Moreover, the relation on the set of slices defined by   $S<S'$  if $\forall x\in S, y\in S'$, $\Diff_{a,b} (x,y) <0$, is a total order on the set of slices.
\end{lemma}

{\it Proof. } All the assertions follow immediately from the first one, which follows from a short computation   (we reproduce that of  \cite[Lemma 3.4]{RS}).  Notice that (writting $N$ for $N^{a,b}$): $\Card (N_L(x)\setminus N_L(y)) - \Card (N_L(y)\setminus N_L(x)) + \Card (N_L(y)\setminus N_L(z)) -\Card (N_L(z)\setminus N_L(y))$ is equal to $\Card(N_L(x)\setminus N_L(z)) - \Card(N_L(z)\setminus N_L(x)) $, and similarly for $N_R$. $\square$

\begin{lemma}[Properties of slices]\label{lem;bloc}
\begin{itemize}
\item[(i)] If $v\in Cyl(a,b)$, then 
 $v$   is at distance at most $2\delta$ from any tight  geodesic segment $[a,b]$. 
\item[(ii)] Let $S$ be a slice of $Cyl(a,b)$, and  $v, v'$ in $S$. Then $d(v,v')\leq
200\delta$.
\item[(iii)]  Let $v$ and $v'$ be in two consecutive slices of $Cyl(a,b)$. Then $|v-v'|\leq 1000 \delta$.
\item[(iv)]  If  $Cyl(a,b) \cap B_{R}(a) = Cyl(a,c) \cap B_{R}(a)$ 
(where $B_R(a)$ is the ball centered at $a$ of radius $R$), 
then any slice of $Cyl(a,b)$ included in
$B_{R-200\delta}(a)$ is a slice of $Cyl(a,c)$.
\end{itemize}
\end{lemma}

{\it Proof. } (Here is a repetition of the proofs of Lemma 2.19-2.21 from \cite{Dis}). For $(i)$ it suffices to see that a point in a cylinder is in a geodesic starting and ending at distance $2\epsilon$ 
from  $[a,b]$, and sufficiently far from its endpoints. 

For $(ii)$,  assume that $d(a,v) \leq d(a,v')$, and $d(v,v')\geq 200\delta$, then the result follows from the relations 
$ N_L^{(a,b)} (v) \subset N_L^{(a,b)} (v') $ (strict inclusion), and $  N_R^{(a,b)} (v') \subset N_R^{(a,b)} (v)$. We now prove the first one (the second one is similar).
The equality is impossible since $v'$ is in one and not the other. Let $w\in [a,b]$ at distance $2\delta$ from $v$, and similarly $w'\in [a,b]$ close to $v'$. Clearly $d(a,w) \leq d(a,w') -196\delta$. If $z\in  N_L^{(a,b)} (v)$, there is $w_z$ on $[a,b]$ at distance $2\delta$ from it. By definition of   $N_L^{(a,b)} (v)  $ it follows that $d(w,w_z) \geq 96\delta$, and $w_z \in [a,w]$. By the triangle inequality, one finds that  $z\in  N_L^{(a,b)} (v')$, what we wanted.

For $(iii)$, assume $d(a,v) < d(a,v')$, and $d(v,v')> 1000 \delta$. 
We can find $w\in [a,b]$ at distance at least $400\delta$ from $v$ and from $v'$ and such that $d(a,v)+200\delta < d(a,w) < d(a,v') -200\delta$. It is easy to check that $w$ is in a slice between that of $v$ and that of $v'$, contradicting that they are in consecutive slices.

For $(iv)$,  let  $S, S' \subset B_{R-200\delta}(a) $ be  slices of respectively  $Cyl(a,b)$, and $Cyl(a,c)$. We claim that if they intersect, they are equal.  Let $v  \in S\cap S'$, and $v'\in S$.  It is enough to check that $\Diff_{a,b}(v,v') = \Diff_{a,c}(v,v')$, because this would implies $v'\in S'$, and $S\subset S'$, and by symmetry, equality. By assumption  $Cyl(a,b) \cap B_{R}(a) = Cyl(a,c) \cap B_{R}(a)$ therefore $N_L^{(a,c)} (v) = N_L^{(a,b)} (v)$, and similarly for $v'$. If $x\in N_R^{(a,c)} (v') \setminus  N_R^{(a,c)} (v)$, it is $100\delta$-close to $v$, and it is then in $Cyl(a,b)$, and in  $N_R^{(a,b)} (v') \setminus  N_R^{(a,b)} (v)$.  By symmetry we also have the reverse inclusion, and  $N_R^{(a,c)} (v') \setminus  N_R^{(a,c)} (v) = N_R^{(a,b)} (v') \setminus  N_R^{(a,b)} (v)$, which ensures that  $\Diff_{a,b}(v,v') = \Diff_{a,c}(v,v')$. $\square$

As a consequence of Theorem \ref{theo;canonicalcylinders} and  Lemma
\ref{lem;bloc} {\it (iv)}, one gets:

\begin{prop}\label{prop;slices}

 We keep the notations of Theorem \ref{theo;canonicalcylinders}, and let $l$ be the constant  given by it. 
  Let $(x,y,z) = (p,
\alpha p, \gamma^{-1}p)$  be a triangle in $\mathcal{CG}(\Sigma)$, such that $\alpha,
\beta, \gamma$ are in $F \cup F^{-1}$, and $\alpha \beta \gamma =1$.

The ordered slice decomposition of the cylinders is as follows.
$$
\begin{array}{lllllllllll}
 Cyl_l(x,y) & = &(S_1, & S_2,  & \dots, & S_k, &\, \mathcal{H}_z,\, & T_m,  & T_{m-1}, & \dots, & T_1)   \\
 Cyl_l(x,z) & = &(S_1, & S_2,  & \dots, & S_k, &\, \mathcal{H}_y,\, & V_p, & V_{p-1}, & \dots, & V_1)   \\
 Cyl_l(y,z) & = &(T_1, & T_2,  & \dots, & T_m, &\, \mathcal{H}_x,\, & V_p, & V_{p-1}, & \dots, & V_1),
\end{array}
$$

 such that $S_1, \dots, S_k, T_1,\dots, T_m$ and  $V_1,\dots, V_p$ are slices and that each $\mathcal{H}_v$,
$(v=x,y,z)$ is a set of at most $10 \psi(n)$ consecutive slices.
The sets  $\mathcal{H}_v$ are called the \emph{holes} of the slice decomposition.

\end{prop}

\section{Purely pseudo-Anosov images of groups}\label{sec;images}

To get Theorem \ref{theo;intro1}, we can follow the approach in \cite{Del}, without major change.  

Let us consider $\varphi:G \to \mathcal{MCG}(\Sigma)$ and $P(G)$ a Van Kampen 2-complex of $G$ (so, $G\simeq \pi_1 (P(G))$ once a base point is chosen). 
 
 $P(G)$ is a simplicial complex with one vertex (and base point), and a certain number of edges $e_1,\dots e_r$, that we identify with elements of $G$, and $T(G)$ triangles (which are the relators of a triangular presentation)\footnote{In principle we could use relations of length $2$ or $3$, but since pseudo-Anosov elements of $\MCG$ are of infinite order, we can assume $G$ to be torsion free, in particular without element of order $2$, and it is not hard then to eliminate relations of length $2$.}. Then we pass to the universal cover
$\widetilde{P(G)}$, where we choose a base point, and representatives $\tilde{e_i}$ of the $e_i$, 
	starting at this point. The vertices of $\widetilde{P(G)}$ are
thus identified with the group $\pi_1(P(G))$, and the map $\varphi$ induces a map
from the vertices of $\widetilde{P(G)}$  to $\MCG(\Sigma)$, and therefore, by
considering the orbit of the base point $p \in \CG(\Sigma)$, to the curve
graph $\CG(\Sigma)$. Denote this map by $\tilde{\varphi}  : \widetilde{P(G)}^{(0)} \to
\CG(\Sigma)$. We then apply  Proposition \ref{prop;slices} to the family  $\varphi(e_1),\dots \varphi(e_r)$, thus providing canonical cylinders for each pair $(p,  \tilde{\varphi}(e_i)p)$ in $\CG(\Sigma)$, 
	and their translates (we will omit the constant $l$ in notations, since it is now fixed until the end). We then extend the map $\tilde{\varphi}$ to the $1$-skeleton of $\widetilde{P(G)}$, by mapping, for all $i$, the edge   $\tilde{e_i}$ onto a path in $Cyl(p,  {\varphi}(e_i)p)$ that successively go through all consecutive slices, and then extend the map $\tilde{\varphi}$ on the translates of $\tilde{e_i}$ equivariantly. This gives an equivariant map    $\tilde{\varphi}  : \widetilde{P(G)}^{(1)} \to
\CG(\Sigma)$.

If $Cyl(p,  {\varphi}(e_i)p)$ has $n$ slices, we choose $n$ distinct points $m_k$, $k=1,\dots, n$, (and we call them ``marked points'') on the edge $\tilde{e_i}$, such that for all $k$, $\tilde{\varphi}(m_k)$ is in the $k$-th slice of $Cyl(p, {\varphi}(e_i)p)$. We complete by translation, so that every edge $\tilde{e}=(v_1,v_2)$ of  $\tilde{S_g}$ has a certain number of points marked on it, that are mapped into the consecutive slices of the cylinder $Cyl(\tilde{\varphi}(v_1),  \tilde{\varphi}(v_2))$. We will use the coincidence of slices to construct tracks in $P(G)$.

	Let us  consider a representative of an orbit of triangular cells in  $\widetilde{P(G)}$. We link each pair of marked points on the edges by a ``blue'' segment,  when the slices in which they are mapped  are equal.  After that, in each triangle where there are unlinked marked points, we add a singular red point in the triangle, and link it with every remaining marked point, by red edges that do not cross any blue one (it is clear that there is a way of choosing the red point so that this is possible). By Proposition \ref{prop;slices}, each singular red point is linked by a red edge to at most $30\psi(T(G))$ marked points.

The union $\mathcal{BR}$ of these segments defines a family of disjoint connected graphs in   $\widetilde{P(G)}$, some graphs with blue and red edges, and some graphs with only blue edges.

	We now extend the map  $\tilde{\varphi}$ on each of these graphs. It suffices to choose the image of each  blue or red segment joining two marked points in a triangle, and then complete by translations. We thus choose any path in  $\CG(\Sigma)$ joining the images of its end points (and similarly for red ones). 	By construction this map  $\tilde{\varphi}  : \widetilde{P(G)}^{(1)}\cup \mathcal{BR} \to
\CG(\Sigma)$ is still equivariant, meaning that $\tilde{\varphi}(gv) = \varphi(g) \tilde{\varphi}(v)$ for all $v\in \widetilde{P(G)}$ and all $g\in G$.  Hence we have:

 \begin{lemma}\label{lem;blueleaf}  The image of a connected completely blue graph of  $\mathcal{BR}$ is contained in a single slice, and thus is finite. 
\end{lemma}

	Since the construction was done $G$-equivariantly, $\mathcal{BR}$ descends to the quotient  $P(G)$ as $\overline{\mathcal{BR}}$ which is the union of disjoint connected graphs, some of them completely blue, some of them containing red edges. Since there are $T(G)$ triangles in $P(G)$, there are at most $T(G)\times 30 \psi(T(G))$ red edges in   $\overline{\mathcal{BR}}$.

\begin{figure}[hbt]
            \begin{center}
              \includegraphics{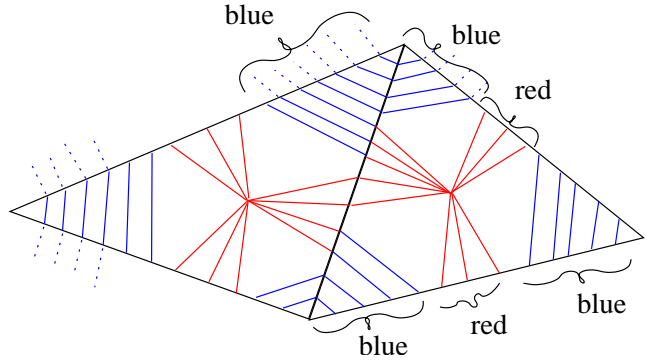}
              \caption{Two triangles in $\widetilde{P(G)}$, with blue and red segments. Consecutive blue segments are mapped into a single slice in $\CG(\Sigma)$, and red ones are mapped on segments of controlled length.  Note that blue segments in a triangle can turn red in another} 
                \label{fig;inPtilde}
           \end{center}
         \end{figure}

To illustrate our construction, let us describe  $\overline{\mathcal{BR}} $ in a triangle $T$ of $P(G)$. 
  The components of $T\smallsetminus   \overline{\mathcal{BR}}$ are either triangles with one blue edge (the spikes around the vertices of $T$, there are three of them, given $T$), or quadrilaterals with two parallel blue edges (there are arbitrarily many of those) or triangles with two red edges, and one vertex being the singular red point of the triangle (there are at most $3\times 10\psi(T(G))$ of them),  or pentagons with two consecutive red edges (containing the red singular point) and the opposite edge blue (there are three of them, given $T$).

\begin{lemma}\label{lem;dC}

To every component
$C$ of $P(G)$ one may associate a graph $dC$ with blue and red edges,
having one or two connected components, so that the following hold:

\begin{itemize}
\item[(i)] In the disjoint union of all graphs  $dC$, over all components $C$,  there are  at most $2 T(G)\times  30\psi(T(G))$ red edges (twice the total number of red edges in  $\overline{\mathcal{BR}}$).

\item[(ii)] Let $C$ be a component, and  $\partial C= \bar{C} \setminus C$. Any loop $l_1$ in $dC$ is homotopic to a loop $l_2$ in $\partial C$ such that $l_1$ has, for each color,  at most twice the number of edges $l_2$ has.    

\item[(iii)] The embedding of a component of $dC$ in $C$ induces an epimorphism on the $\pi_1$  except possibly if $dC$ has one component, in which case $\pi_1(dC)$  maps   
in $\pi_1(C)$  either surjectively, or in a subgroup of index $2$.

\item[(iv)] The group  $\pi_1(C)$ is free.

\end{itemize}
\end{lemma}

{\it Proof. } Let us consider $C$ a connected component of  $P(G) \smallsetminus \overline{\mathcal{BR}}$. In $C$ we choose $dC$ to be the boundary of a small tubular neighborhood of $\partial C = \bar{C} \setminus C$. This is a graph, and for every triangle $T$ of $P(G)$,     and any component $C_0$ of $T\cap C$,   $dC\cap C_0$ has same number of components as $\partial C \cap \overline{C_0}$ has (that is one or two).  We color each of them by the color of the neighboring component of   $\partial C \cap  \overline{C_0}$, thus ensuring $(i)$, since every  edge of    $\overline{\mathcal{BR}}$  locally separates $P(G)$ in two.

  When homotoping $dC$ to $\partial C$ in the relevant tubular neighborhood of $\partial C$, we send edges of a given color on edges of $\partial C$ of the same color, and, again because  every  edge of    $\overline{\mathcal{BR}}$  locally separates $P(G)$ in two, at most two edges of $dC$ on the same edge of $\partial C$, hence $(ii)$.

Let now $p$ be a base point in $dC$, and $\ell$ a loop in $C$ starting at $p$. Let us denote by $T_1, \dots, T_k=T_1$ the consecutive triangles in which $\ell$ enters. For $i\leq k-1$, we can inductively homotope $\ell$ so that it stays in the same component of $dC \cap T_i$. Then in $T_k$, either $\ell$ enters in the component of $p$, or it enters in the other component (if any). In the first case, $\ell$ is homotopic to a loop of $dC$, and in the second case, $\ell^2$ is homotopic to a loop of $dC$, moreover in this second case, $dC$ (globally) has only one connected component.

Thus, in the first case, the inclusion induces an epimorphism on the fundamental groups. In the second case, we need to show that, if   $\pi_1(dC)$ is not surjective in $\pi_1(C)$, then it is contained in a subgroup of index $2$.  The image of  $\pi_1(dC)$ contains the subgroup $S$ generated by all the squares of $\pi_1(C)$. This subgroup $S$ is normal, since a conjugate of a product of squares is a product of squares, and the quotient $\pi_1(C) / S$ is finitely generated (as $\pi_1(C)$ is), and has all its elements of order $2$. Hence it is abelian, hence finite, and it is isomorphic to  $(\mathbb{Z}/2\mathbb{Z})^n$ for some $n$. Now, if  $\pi_1(dC)$ is not surjective in $\pi_1(C)$, since it contains $S$, it maps on $\pi_1(C) / S$ on some proper subgroup. Since the quotient is abelian, one can quotient by this proper subgroup, and this gives a certain  $(\mathbb{Z}/2\mathbb{Z})^k$, $k\leq n$. Hence there is a surjective map on $\mathbb{Z}/2\mathbb{Z}$ that contains $\pi_1(dC)$ in its kernel, what we wanted.
 This establishes point $(iii)$.  

To obtain $(iv)$, it suffices to notice that $C$ is homotopically equivalent to $C' = C \setminus \{C_i\}$, where the $C_i$ are the components of $C\cap T$, for some $T$,  that are triangles.  Now $C'$ is a union of  quadrilaterals and pentagons, glued together on two opposite sides. Consider the graph with one edge in each such quadrilateral or pentagon,  joining the midpoints of the opposite sides on which the gluing is done.  It is easily checked that  $C'$ is then homeomorphic to a fibration of an open interval $(-1,1)$ 
on this graph, which makes its fundamental group free, hence the fourth point. $\square \smallskip$

\begin{lemma}

Let $X$ be the bipartite graph such that the vertices of one color (white) are the components of $\overline{\mathcal{BR}}$, the vertices of the other color (black) are the components  of $P(G) \setminus \overline{\mathcal{BR}}$, and that the edges realize the adjacency relation.

Then $G$ has a natural structure of fundamental group of a finite graph of groups such that the underlying graph is $X$, the vertex groups are the fundamental groups of the relevant components, and the edge groups
are  the fundamental subgroups of the relevant intersections, 
with identifications, through edges, induced by the adjacency in $P(G)$.

\end{lemma}

{\it Remark.} The graph $X$ is endowed with groups for each vertex and each edge, and attaching maps from edge groups to adjacent vertex group. Here, one should note that the attaching  maps need not be injective, and the vertex and edge groups need not embed in $G$. This is not important if one is interested in finding a presentation of $G$ using this construction. If one wants to find a graph of groups with all maps injective, one should take for vertex groups, instead of the fundamental groups of the components,  their \emph{images} in $G$.

{\it Proof. } This is an application of the Van Kampen theorem, for our decomposition of $P(G)$. $\square \smallskip$

Let us now describe a particular generating set for the vertex and edge groups of this graph of groups.

Let $c$ be a  component of $\overline{\mathcal{BR}}$ (it is a red-and-blue graph in $P(G)$).  
We make a careful choice of generators (this construction actually works for any red-and-blue graph). 

Let $\mathcal{B}$ be the maximal blue subgraph of $c$ (not necessarily connected). Let $b_1, \dots b_k$ be a minimal collection of blue edges (possibly empty) so that $\mathcal{B} \smallsetminus \{ b_1, \dots , b_k\}$ is a forest. By minimality of $k$,  $\mathcal{B} \smallsetminus \{ b_1, \dots , b_k\}$ has same number of connected components as $\mathcal{B}$,  hence  $c \smallsetminus \{ b_1, \dots , b_k\}$  is connected. Let now $r_1, \dots r_n$ be a minimal collection (possibly empty) of red edges so that  $c \smallsetminus \{ b_1, \dots , b_k, r_1, \dots r_n\}$ is a tree. 
Such a collection exists, since if a red-and-blue graph is not a tree, it has a cycle, which cannot consist only of blue edges, if the maximal blue subgraph is a forest. Hence this red edge can be removed, and the new graph is still connected (and has less edges).

Note that  $c \smallsetminus \{ b_1, \dots , b_k, r_1, \dots r_n\}$ is in fact a maximal subtree of $c$, for, if we put back one red edge, it is not a tree, and if we put back one blue edge, some blue component is not a tree. 

Therefore, if $*$ is a base point in $c$, $\pi_1(c,*)$ has a natural isomorphism with the free group on $ \{ b_1, \dots , b_k, r_1, \dots r_n\}$. 
We call the $r_i$ the {\it red generators}, and the $b_j$ the {\it blue} ones.

Let now $C$ be  a component  of $P(G) \setminus \overline{\mathcal{BR}}$. Each component of $dC$ (see Lemma \ref{lem;dC}) is a red-and-blue graph, so the construction above can be performed, thus providing a system of generators of $dC$.

\begin{figure}[hbt]
            \begin{center}
              \includegraphics{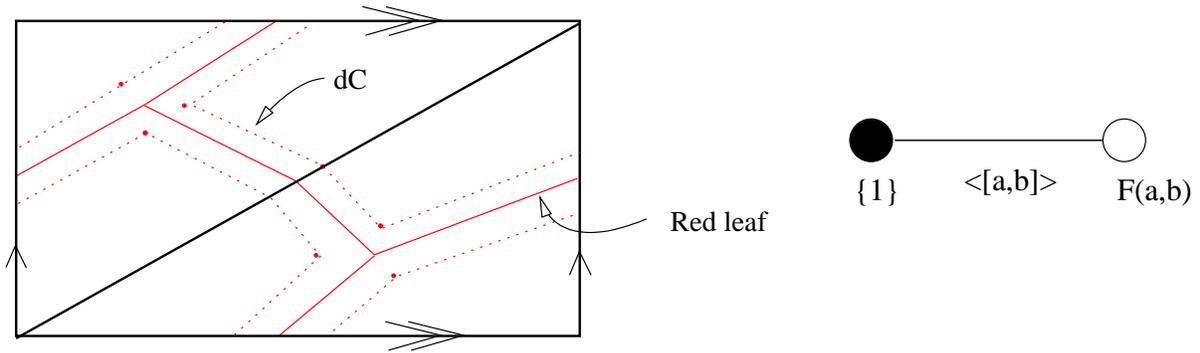}
              \caption{In the torus (left), with only one red leaf (surrounded by  its tubular neighborhood, defining $dC$) the graph $X$ has a white vertex (corresponding to the red graph) and a black vertex (the $2$-cell attached to it)}
                \label{fig;torus}
           \end{center}
         \end{figure}

\begin{lemma} \label{lem;blue} 

If the image of $G$ in $\MCG$ is purely pseudo-Anosov, 
then a loop defining a  blue generator (in the graph $dC$ of a component $C$, or in a component of  $\overline{\mathcal{BR}}$) is trivial in $G$.

\end{lemma}

{\it Proof. } Let  $\gamma$ be a closed blue curve in $\overline{\mathcal{BR}}$, note that $\varphi(\gamma)$ is defined up to conjugacy.  By  Lemma \ref{lem;blueleaf}  $\varphi(\gamma)$ has a finite orbit in $\CG(\Sigma)$ (contained in the slice associated to the blue graph we started from). Thus  $\varphi(\gamma)$ is not pseudo-Anosov, and by assumption this means that $\gamma$ is trivial in $G$.

It remains to check that a loop defining a blue generator is freely homotopic to a blue curve in $P(G)$.  Indeed, by definition of the generator, the loop consists of a path $p$ going from the base point to a vertex of a blue edge, then this blue edge $e$ and a path $q$, with the property that $p$ and $q$ are in a same maximal subtree $T$ of the graph, containing a maximal blue subtree $T_b(e)$ of the blue graph $B(e)$ containing $e$. We claim that $p$ and $\bar{q}$ (with reverse orientation) enter $T_b(e)$  on the same point. If it was not the case, a path $p_b$ in $T_b(e)$ between the two entering points would give rise to a loop $(p) (p_b) (q)$ in $T$, contradicting that $T$ is a tree. This implies also that from the base point to this entering point, the paths $p$ and $\bar{q}$ are equal. Hence, the loop is freely homotopic to a loop contained in the graph $G(e)$ which is completely blue. $\square \smallskip$

In the following $b_1(G, \mathbb{Z}/2\mathbb{Z})$ is the first Betti number of $G$ over $\mathbb{Z}/2\mathbb{Z}$. Let $N_0$ be the number of red edges in  $\overline{\mathcal{BR}}$ in $P(G)$. Note that $N_0\leq  T(G)\times 30 \psi(T(G))  $ (see remark after Lemma \ref{lem;blueleaf}).

\begin{lemma}\label{lem;graph}  $G$ has a presentation as a bipartite graph of groups $X'$  (with black and white vertices as above), that satisfies  the following properties.
 
\begin{itemize}
\item There are at most $ T(G)$ white vertices, whose groups are free with at most $N_0$ generators each.  
\item There are at most  $2 N_0 $ black   vertices  
whose groups are free with at most    $2 N_0 $ generators each. 
\item For each edge, and each adjacent vertex,   the   
  attaching  map from the edge group into the vertex group sends each generator to a product of at most $2N_0$ of the given generators.  
\item For each edge  the 
 attaching map from the edge group into the adjacent black vertex group 
 is either surjective or has its image in a subgroup of index $2$.

\item  There are  at most  $b_1(G,\mathbb{Z}/2\mathbb{Z})$ edges for which the 
 attaching map into the adjacent black vertex group is not surjective.

\end{itemize}

\end{lemma}

{\it Proof. } The graph of groups $X'$ is not necessarily $X$ defined above: first we simply remove every vertex of $X$ with no red generator (and their adjacent edges), since by Lemma \ref{lem;blue} their groups are trivial. Then, if an edge $e$ 
 of $X$ has only blue generators (hence with trivial group), it is separating (otherwise $G$ would have a cyclic free factor, hence  several ends), and one of the components of $X \setminus e$ has trivial group (otherwise $G$ has several ends).  By removing this subgraph, we can assume that no edge has only blue generators without changing the fundamental group. Thus we get the graph of groups $X'$, and  only now it is possible to bound the number of white vertices and black vertices,  respectively by the total number of components of $\overline{\mathcal{BR}}$ with a red edge (less than the total number of red singular points, $T(G)$), and by twice the number of red edges in   $\overline{\mathcal{BR}}$ (each component of $P(G) \setminus \overline{\mathcal{BR}}$ associated to a vertex of  $X'$ is adjacent to a certain red edge, and only two of them can be adjacent to the same red edge).  

The generators of each vertex or edge groups are chosen to be the red ones in the construction above (the blue ones being all trivial by Lemma \ref{lem;blue}). The obtained graph satisfies then the two first  points, by construction. 

Given an edge, it corresponds to an adjacency of a component $C$ of $P(G) \setminus \overline{\mathcal{BR}}$ (its black vertex) and a component of $\overline{\mathcal{BR}}$ (its white vertex). Hence it corresponds to a component of $dC$.

By Lemma \ref{lem;dC}$(ii)$, each  loop in a component of $dC$  is homotopic to a loop in the relevant component of  $\overline{\mathcal{BR}}$ containing at most twice the number of edges of each color. If the loop in $dC$ is simple, it is homotopic to a loop  of  $\overline{\mathcal{BR}}$ passing at most twice through each edge it contains. Thus it is homotopic to  the product of $2N_0$ red generators of the relevant component of  $\overline{\mathcal{BR}}$, which proves the third point for attaching maps into white vertex groups. For attaching maps into black vertex groups,  the bound is obtained similarly, replacing the component of $\overline{\mathcal{BR}}$ by the graph in $C$ to which $C$ is homotopically equivalent.

The fourth point follows  by Lemma \ref{lem;dC}$(iii)$.

It remains to bound the number of edges whose group is in a subgroup of index $2$ in the adjacent black vertex group: each of them gives rise to a morphism of $G$ onto $\mathbb{Z}/2\mathbb{Z}$ (by sending  the index $2$ subgroup containing the edge group on $0$, and  its non-trivial coset in the white vertex on $1$). One cannot have more than  $b_1(G,\mathbb{Z}/2\mathbb{Z})$ distinct such morphisms. $\square \smallskip$

In the graph $X'$ we choose  a base white vertex $v_0$, and in $P(G)$, a base point $p_0$ in the component of $\overline{\mathcal{BR}}$ of $v_0$. 
 For each black vertex $w$ of $T'$ (which has valence at most $2$), we choose a path between its two adjacent components of $\overline{\mathcal{BR}}$ as follows.   By adjacency, there is a triangle of the initial triangulation of $P(G)$ in which these components have adjacent segments. We choose a segment $s_w$ between two such points in that triangle, so that $s_w$ does not intersect any other component. 

 We thus can choose a simple path from $p_0$ to any component of $\overline{\mathcal{BR}}$ by a sequence of simple paths in  components, and segments  $s_w$ for black vertices $w$ of $X'$. Once chosen a maximal subtree $T'$ in $X'$, this gives a choice of one simple path from $p_0$ to any component, hence a choice of well defined conjugacy classes of any vertex group of $X'$ in $G$. 

Thus, given a component $c$ of   $\overline{\mathcal{BR}}$, its base point $p$ is chosen to be the first point  of the component met by the former path, and the red generators of $\pi_1(c,p)$ are seen as elements of $G$. 

Finally, we choose in the universal cover $\widetilde{P(G)}$ a pre-image $\tilde{p}_0$ of $p_0$, and from it, cross sections of each of the chosen path from $p_0$ to a component of $\overline{\mathcal{BR}}$. If $p$ is a base point of a component, this gives $\tilde{p}$ a particular pre-image of $p$ in $ \widetilde{P(G)}$.

The next lemma should be compared with \cite[Lemma III.4 ]{Del}.
\begin{lemma}\label{lem;del} 

Let $c$ be a  component of  $\overline{\mathcal{BR}}$ that contains a red edge, and $p$ the base point in $c$. 
 Let $v \in \CG(\Sigma)$ be  $v= \tilde{\varphi} (\tilde{p})$ with the notation just introduced.
 
Then 
 the image under $\varphi$ of each of the red generators of $\pi_1(c)$  translates $v$  
by at most $T(G)( 20\psi(T(G))\times 1000\delta  +200\delta )$  in  $\CG(\Sigma)$.

Moreover, if $c'$ is another component that is adjacent to $c$ in $T'$ (their white vertices in $T'$ are at distance $2$), then the vertex $v'$ for $c'$ is 
at distance at most $1000\delta + T(G)( 20\psi(T(G))\times 1000\delta  +200\delta )$ from $v$, in $\CG(\Sigma)$.
 
\end{lemma}

{\it Proof. } Let us choose a connected cross section  $\tilde{c}$ of $c$ in $\widetilde{P(G)}$, from the point $\tilde{p}$. 

It is a graph with blue and red edges, and $\tilde{\varphi}$ maps it in  $\CG(\Sigma)$. Any red generator of $\pi_1(c)$ move the base point of $ \tilde{c}$
to a point of $\tilde{c}$. 
 The equivariance of $\varphi$ implies that the image in $\mathcal{MCG}(\Sigma)$ of the generators of $\pi_1(c)$ move  $v= \tilde{\varphi}(\tilde{p})$ (which we choose as the vertex $v$ of the statement)
to a point of $\tilde{\varphi}(\tilde{c})$.

  Each blue segment of it is mapped in a finite subgraph (a slice in fact), with universally bounded diameter (at most $200\delta$, by Lemma \ref{lem;bloc} {\it (ii)}).

Each pair of red edges around a singular point (corresponding to the center of the slice decomposition of a triangle) is mapped on a path between to points at uniformly bounded distance (at most $20\psi(T(G))\times 1000\delta$, by Proposition \ref{prop;slices} and Lemma \ref{lem;bloc} {\it (iii)}).

Now in any segment of  $\tilde{c}$, there are at most $T(G)$ different singular points, since it is a cross section of a graph on $P(G)$. Thus there are at most $T(G)$ pairs of red edges as above. Therefore, the extremal points of a segment in   $\tilde{c}$ are mapped by $\tilde{\varphi}$ to points at distance at most $T(G)( 20\psi(T(G))\times 1000\delta  +200\delta )$ from each other. This proves the first claim.

Now if $c'$ is a component  adjacent to $c$ in $T'$, its base point is mapped in a slice that is adjacent to a slice of a point of $\tilde{\varphi}(\tilde{c})$, thus at distance at most $1000\delta$ from this point, by  Lemma \ref{lem;bloc} {\it (iii)}). With the former estimate on the diameter of $\tilde{\varphi}(\tilde{c})$, this gives the required bound.
$\square \smallskip$

We can now prove  Theorem  \ref{theo;intro1}.

{\it Proof. } Let $D_0 =  T(G) \times(20\psi(T(G))\times 1000\delta + 200\delta)$ 
and $\Diam(X')$ be the diameter of the graph $X'$; these two constants depend only on $G$ and $\Sigma$.
The generating set is that given by the graph of group $X'$ of  Lemma \ref{lem;graph}, taking red generators for every vertex and edge groups. From  Lemma \ref{lem;graph}, we can write a presentation over this generating set: the relations are the words $s^{-1}f(s)$ where $s$ runs over  the generating sets of the edge groups, and $f(s)$ is the image of $s$ under the corresponding attaching
map (and note $f(s)$ is equal to a product of at most $2N_0$ generators
of the range of this map).

Hence we have a bound on the complexity of the presentation over this generating set. From  Lemma \ref{lem;del} we get a subset $\Delta$ of $\CG(\Sigma)$ of diameter bounded by $\Diam(X') \times (D_0+1000\delta)$ such that for all generator $r$ in our family, there is $v\in \Delta$ such that $d(v,\varphi(r) v)$  is universally bounded (by $D_0$ if the generator is in a vertex group, and by $\Diam(X') \times (D_0+1000\delta)$  if the generator is a stable letter of the graph $X'$ with the maximal subtree $T'$). The triangle inequality easily implies that for any point $v$ in $\Delta$, and any generator $r$ in our family, the displacement is bounded by $2\Diam(\Delta) +  \Diam(X') \times (D_0+1000\delta)$ which depends only on $G$ and $\Sigma$. $\square$

\medskip

{\sc Francois Dahmani, Institut de Math\'ematiques de Toulouse, Universit\'e Paul Sabatier (Toulouse III)
31062 Toulouse, Cedex 9,    France.}

{\tt e-mail: } { \it francois.dahmani@math.univ-toulouse.fr}

\medskip

{\sc Koji Fujiwara,
Graduate School of Information Science, Tohoku University,
Sendai, 980-8579, Japan.}

{\tt e-mail: } {\it fujiwara@math.is.tohoku.ac.jp}
\end{document}